\documentclass[10pt,a4paper]{article}
\usepackage{graphics, graphicx}     

\usepackage{cite}
\usepackage{amsmath}
\usepackage{amssymb}
\usepackage{amsthm}
\theoremstyle{plain}
\newtheorem{Theorem}{Theorem}[section] %
\newtheorem{Lemma}{Lemma}[section]

\theoremstyle{definition}
\newtheorem{Remark}{Remark}[section]

\theoremstyle{definition}

\newtheorem{Problem}{Problem}[section]

\newenvironment{Proof} 
{\par\noindent{\it Proof of}} 
{\hfill$\vspace{5mm}\scriptstyle\blacksquare$} 

\numberwithin{equation}{section} 
\numberwithin{figure}{section} 
\numberwithin{table}{section} 

\begin{document}

\setcounter{page}{1}

\markboth{M.I. Isaev, R.G. Novikov}{Effectivized H\"older-logarithmic stability estimates}

\title{Effectivized H\"older-logarithmic stability estimates for the  Gel'fand inverse problem}
\date{}
\author{ { M.I. Isaev and R.G. Novikov}}

\maketitle
\begin{abstract}
	We give effectivized H\"older-logarithmic energy and regularity dependent stability estimates
	for the Gel'fand inverse boundary value problem in dimension $d=3$. This effectivization includes 
	explicit dependance of the estimates on coefficient norms and related parameters. Our 
	new estimates are given in $L^2$ and $L^\infty$ norms for the coefficient difference and related stability efficiently 
	increases with increasing energy and/or coefficient difference regularity. Comparisons with preceeding results are given.
\end{abstract}

\section{Introduction and main results}
 We consider the equation
\begin{equation}\label{eq} 
	-\Delta \psi  + v(x)\psi = E\psi, \ \  x \in  D \subset \mathbb{R}^3,
\end{equation}
where
\begin{equation}\label{eq_c}
	\begin{aligned}
			D  \text{ is an open bounded domain in $\mathbb{R}^3$}, \ \ \partial D \in C^2,
 	\end{aligned}
\end{equation}
\begin{equation}\label{eq_c1}
	 v\in L^{\infty}(D).
\end{equation}
Equation (\ref{eq}) can be regarded as the stationary  Schr\"odinger equation of quantum mechanics
 at fixed energy  $E$. Equation (\ref{eq}) at fixed $E$ arises also in acoustics and electrodynamics.

As in Section 5 of Gel'fand's work \cite{Gelfand1954} we consider an operator establishing a relationship between $\psi$ 
and $\partial \psi / \partial \nu$ on $\partial D$ for all sufficiently regular solutions $\psi$ of equation (\ref{eq}) in $\bar{D} = D \cup \partial D$
at fixed $E$, where $\nu$ is the outward normal to $\partial D$.
As in \cite{Novikov1988}, \cite{IN2012} (for example) we represent such an operator as the Dirichlet-to-Neumann map  $\hat{\Phi}(E)$ defined by the relation
\begin{equation}\label{def_phi}
	\hat{\Phi}(E) (\psi|_{\partial D}) = \frac{\partial \psi}{\partial \nu}|_{\partial D},
\end{equation}
where we assume also that
\begin{equation}\label{correct}
	\text{$E$ is not a Dirichlet eigenvalue for  operator $-\Delta + v$ in $D$.} 
\end{equation}
The map $\hat{\Phi} = \hat{\Phi}(E)$  can be regarded
as all possible boundary measurements for the physical model described by equation \eqref{eq} at fixed energy $E$
under assumption \eqref{correct}.

We consider the following inverse boundary value problem for equation (\ref{eq}): 

\begin{Problem}
 Given $\hat{\Phi}$ for some fixed $E$, find $v$. 
\end{Problem}

This problem is known as the Gel'fand inverse boundary value problem for the Schr\"odinger equation at fixed energy $E$ in three dimensions (see \cite{Gelfand1954}, \cite{Novikov1988}). For $E=0$ this problem can be regarded also as a generalization of the Calder\'on problem of the electrical impedance tomography in three dimensions (see \cite{Calderon1980}, \cite{Novikov1988}). 
Problem 1.1 can be also considered as an example of ill-posed problem; 
see  \cite{BK2012}, \cite{LR1986} for an introduction to this theory.

 Let, for real $m\geq 0$,
\begin{equation}\label{def_Hm}
	\begin{aligned}
	H^{m}(\mathbb{R}^3) = \left\{w \in L^2(\mathbb{R}^3):\ \mathcal{F}^{-1} (1+|\xi|^2)^{\frac{m}{2}} \mathcal{F} w \in  L^2(\mathbb{R}^3) \right\},
\\
	||w||_{H^{m}(\mathbb{R}^3)} = \left\|\mathcal{F}^{-1} (1+|\xi|^2)^{\frac{m}{2}} \mathcal{F} w \right\|_{L^2(\mathbb{R}^3)},
\end{aligned}
\end{equation}
where $\mathcal{F}$  denote the Fourier transform
\begin{equation*}
\mathcal{F}w (\xi) = \frac{1}{(2\pi)^3}\int\limits_{\mathbb{R}^3} e^{i\xi x} w(x) dx, \ \ \ \xi\in \mathbb{R}^3.
\end{equation*}

In addition, for real $m\geq 0$, we consider the spaces $W^{m}(\mathbb{R}^3)$ defined by 
\begin{equation}\label{def_Wm}
	\begin{aligned}
	W^{m}(\mathbb{R}^3) = \left\{w \in L^1(\mathbb{R}^3):\ (1+|\xi|^2)^{\frac{m}{2}} \mathcal{F} w \in  L^\infty(\mathbb{R}^3) \right\},
\\
	||w||_{W^{m}(\mathbb{R}^3)} = \left\| (1+|\xi|^2)^{\frac{m}{2}} \mathcal{F} w \right\|_{L^\infty(\mathbb{R}^3)}.
\end{aligned}
\end{equation}
We note that for integer $m$ the space $W^{m}(\mathbb{R}^3)$ contains 
the standard Sobolev space $W^{m,1}(\mathbb{R}^3)$  of $m$-times smooth functions in $L^1$ on $\mathbb{R}^3$.


In the present work we obtain, in particular, the following theorems:

\begin{Theorem}\label{Theorem1} 
Suppose that $D$ satisfies (\ref{eq_c}) and
 $v_1, v_2$ satisfy  (\ref{eq_c1}), (\ref{correct}) for some real $E$.
 Suppose also that:  $||v_{j}||_{L^\infty(D)} \leq N$ for some $N>0$, $j = 1,2$; 
${\rm supp} (v_2 - v_1) \subset D$, $v_2-v_1 \in H^m(\mathbb{R}^3)$, $\|v_2 - v_1\|_{H^m(\mathbb{R}^3)} 
\leq N_{H^m}$ 
for some $m>0$ and $N_{H^m}>0$. Let  
\begin{equation}\label{Theorem_delta}
	\delta = ||\hat{\Phi}_{2}(E) - \hat{\Phi}_{1}(E) ||_{L^\infty(\partial D)\rightarrow L^\infty(\partial D)},
\end{equation}
where $\hat{\Phi}_{1}(E)$, $\hat{\Phi}_{2}(E)$ denote 
	the Dirichlet-to-Neumann maps for
	$v_1$, $v_2$, respectively.
Then, there exist some positive constants $A,B,\alpha,\beta$ depending on $D$ only such that 
\begin{equation}\label{eq_Theorem1}
	\begin{aligned}
	||v_2 - v_1||_{L^2 (D)} &\leq  
	A\left( \alpha E + \beta (1-\tau)^2\left(\ln\left(3+\delta ^{-1}\right)\right)^2 \right)^{\frac{1}{2}} \delta^\tau + \\ 
	+
	 B  &\left( 1+N\right)^{\frac{4m}{3}} N_{H^m}  \left( \alpha E + \beta(1-\tau)^2 \left(\ln\left(3+\delta ^{-1}\right)\right)^2 \right)^{-\frac{m}{3}} 
	\end{aligned}
\end{equation}
for any $\tau \in(0,1]$ and  $E\geq 0$.

 Besides, estimate (\ref{eq_Theorem1})
is also fulfilled for  any $\tau \in(0,1)$ and  $E<0$ under the following additional condition:
\begin{equation}
	\alpha E + \beta (1-\tau)^2 \left(\ln\left(3+\delta ^{-1}\right)\right)^2  > 0.
\end{equation}
\end{Theorem}

\begin{Theorem}\label{Theorem2} 
Suppose that $D$ satisfies (\ref{eq_c}) and
 $v_1, v_2$ satisfy  (\ref{eq_c1}), (\ref{correct}) for some real $E$.
 Suppose also that:  $||v_{j}||_{L^\infty(D)} \leq N$ for some $N>0$, $j = 1,2$; 
${\rm supp} (v_2 - v_1) \subset D$, $v_2-v_1 \in W^m(\mathbb{R}^3)$, $\|v_2 - v_1\|_{W^m(\mathbb{R}^3)} \leq N_{W^m}$ 
for some $m>3$ and $N_{W^m}>0$.
Let  $\delta$ be defined by (\ref{Theorem_delta}).
Then, there exist some positive  constants $\tilde{A},\tilde{B},\tilde{\alpha},\tilde{\beta}$ depending on $D$ only such that 
\begin{equation}\label{eq_Theorem2}
	\begin{aligned}
	||v_2 - v_1&||_{L^\infty (D)} \leq  
	\tilde{A}\left( \tilde{\alpha} E + \tilde{\beta}(1-\tau)^2 \left(\ln\left(3+\delta ^{-1}\right)\right)^2 \right)^{\frac{1}{2}} \delta^\tau + \\ 
	+
	 \tilde{B}  &\frac{\left( 1+N\right)^{\frac{2(m-3)}{3}} N_{W^m}}{m-3}  \left( \tilde{\alpha} E + \tilde{\beta}(1-\tau)^2 \left(\ln\left(3+\delta ^{-1}\right)\right)^2 \right)^{-\frac{m-3}{6}} 
	\end{aligned}
\end{equation}
for any $\tau \in(0,1]$ and  $E\geq 0$.

Besides, estimate (\ref{eq_Theorem2}) is also fulfilled for any $\tau \in(0,1)$ and $E<0$ 
under the following additional condition 
\begin{equation}
	\tilde{\alpha} E + \tilde{\beta} (1-\tau)^2 \left(\ln\left(3+\delta ^{-1}\right)\right)^2  > 0.
\end{equation}
\end{Theorem}

Theorems \ref{Theorem1} and \ref{Theorem2} are proved in Sections 3 and 4, respectively. These proofs are based on Lemmas 
\ref{Lemma_2.1}, \ref{Lemma_2.2} and \ref{Lemma_2.3} given in Section 2. Then these
proofs are based on the intermediate estimates \eqref{L2_rho}, \eqref{L0_rho} which may be of 
independent interest.


\begin{Remark}
	The estimates of Theorem \ref{Theorem2} can be regarded as a significant effectivization 
	of the following estimates of \cite{IN2012} for the three-dimensional case: 
\begin{equation}\label{eq_t1}
	\begin{aligned}
	||v_2 - v_1||_{L^\infty(D)} \leq 
	C_1(N_m,D,m,E)\left(\ln\left(3+\delta ^{-1}\right)\right)^{-s_1} 
	\end{aligned}
\end{equation}
 for $E \in \mathbb{R}$;  
\begin{equation}\label{eq_t2}
	\begin{aligned}
	||v_2 - v_1||_{L^\infty(D)} &\leq C_2(N_m,D,m,\tau)(1+\sqrt{E}) \delta^\tau + \\&+
	C_3(N_m,D,m,\tau) (1+\sqrt{E})^{s- s_1}\left(\ln\left(3+\delta ^{-1}\right)\right)^{-s}
	\end{aligned}
\end{equation}
 for $E\geq 0$, $\tau \in(0,1)$ and any $s \in [0, s_1]$. Here $\delta$ is defined by (\ref{Theorem_delta}) and $s_1 = (m-3)/3$.

	In addition, estimates \eqref{eq_t1} and \eqref{eq_t2} were obtained in \cite{IN2012} 
	under the assumptions that:  $D$ satisfies (\ref{eq_c}),
 $v_j$ satisfies  \eqref{eq_c1}, (\ref{correct}),  
 ${\rm supp}\, v_j   \subset D$, $v_j \in W^{m,1}(\mathbb{R}^3)$, $\|v_j\|_{W^{m,1}(\mathbb{R}^3)} \leq N_m$, $j=1,2$, 
for some integer $m>3$ and $N_{m}>0$.

Actually, Theorem \ref{Theorem2} was obtained in the framework of finding the dependance of $C_1$, $C_2$, $C_3$ 
of \eqref{eq_t1}, \eqref{eq_t2} on $N_m$, $m$ and $\tau$. One can see that the estimates of Theorem \ref{Theorem2} depend explicitely on coefficient norms 
$N$, $N_{W^m}$ and 
 parameteres  $m$, $\tau$ 
and imply  \eqref{eq_t1}, \eqref{eq_t2} 
with some $C_1$, $C_2$, $C_3$ explicitely dependent on $N_m$, $m$, $\tau$ as a corollary. 
Besides,  in Theorem \ref{Theorem2}
 we do not assume that each of potentials $v_1$, $v_2$ is $m$-times differentiable  and is supported in $D$
 (in a similar way with Theorem 2.1 of \cite{Novikov2013}).
 
 By the way we would like to note also that even for $E=0$ the reduction of H\"older-logarithmic stability estimates like \eqref{eq_Theorem1}, \eqref{eq_Theorem2}
 to pure logarithmic estimates like \eqref{eq_t1} is not optimal for large $m$ because of the following asymptotic formula:
 \begin{equation*} 
 		\sup\limits_{\delta \in (0,1]} \left|\frac{\delta^\tau}{ \left(\ln(3+\delta^{-1}) \right)^{-\mu}}\right| 
 		= O\left( \left(\frac{\mu}{\tau}\right)^\mu e^{-\frac{\mu}{\tau}} \right) \ \text{ as } \mu \rightarrow +\infty.
 \end{equation*}
 In particular, even for $E=0$ the H\"older-logarithmic estimates \eqref{eq_Theorem1}, \eqref{eq_Theorem2}
 are much more informative than their possible pure logarithmic reductions.
\end{Remark}

\begin{Remark}
	Theorem \ref{Theorem1} was obtained as an extention of Theorem  \ref{Theorem2}  to the $L^2$-norm case. In addition,
	it is important to note that the second (``logarithmic'') term
	of the right-hand side of \eqref{eq_Theorem1} is considerably better than the analogous term 
	of \eqref{eq_Theorem2}. In particular,
	\begin{equation*}
		\begin{aligned}
		R = O\left( 
		\left(\ln\left(3+\delta ^{-1}\right)\right)^{-\frac{2m}{3}}
		\right)  \ \ \ \text{ for } \delta \rightarrow 0,\\
		R = O\left(E^{-\frac{m}{3}}\right) \ \ \ \text{ for } E \rightarrow +\infty,
		\end{aligned}
	\end{equation*}
	whereas
	\begin{equation*}
		\begin{aligned}
		\tilde{R} = O\left( 
		\left(\ln\left(3+\delta ^{-1}\right)\right)^{-\frac{m-3}{3}}
		\right)  \ \ \ \text{ for } \delta \rightarrow 0,\\
		\tilde{R} = O\left(E^{-\frac{m-3}{6}}\right) \ \ \ \text{ for } E \rightarrow +\infty,
		\end{aligned}
	\end{equation*}
	where $R$ and $\tilde{R}$ denote the second (``logarithmic'') terms of the right-hand sides of \eqref{eq_Theorem1} and
	\eqref{eq_Theorem2}, respectively.
\end{Remark}


\begin{Remark}
	The estimates of Theorem \ref{Theorem1} should be compared also with the following
	estimate of \cite{INUW2014} for the three-dimensional case:
	\begin{equation}\label{eq_INUW}
		\|v_2-v_1\|_{H^{-m}(\mathbb{R}^3)} \leq C
		\left(E^2 \delta + \left(\sqrt{E} + \ln \delta^{-1}\right)^{-(2m-3)}\right),
	\end{equation}
	where $C = C\left(N_m,D, \text{supp}(v_2-v_1),m \right)>0$, $\|v_j\|_{H^m(D)} \leq N_m$ ($j=1,2$), 
	supp$\,(v_2-v_1) \subset D$, $m> 3/2$, $\delta$ is the distance between the boundary measurements (Cauchy data) for 
	$v_1$, $v_2$ and is, roughly speaking, similar to $\delta$ of \eqref{Theorem_delta} and where $\delta \leq 1/e$.
	
	A principal advantage of \eqref{eq_Theorem1} in comparison with \eqref{eq_INUW} consists
	in estimation $v_2-v_1$ in the $L^2$-norm instead of the $H^{-m}$-norm. Besides,
	estimate \eqref{eq_Theorem1}  depends explicitely  on coefficient norms 
$N$, $N_{W^m}$ and 
 parameteres  $m$, $\tau$ in contrast with  \eqref{eq_INUW}. 
 In addition, in \eqref{eq_Theorem1} we do not assume that each of $v_1$, $v_2$ belongs to $H^m$.
\end{Remark}
\begin{Remark}
In the literature on Problem 1.1 estimates of the form \eqref{eq_t1} are known as global logarithmic stability estimates.
The history of these estimates goes back to \cite{Alessandrini1988} for the case when $s_1 \leq 1$ and 
to \cite{Novikov2011} for the case when $s_1>1$.

In addition,  estimates of the form \eqref{eq_Theorem1}, \eqref{eq_Theorem2}, \eqref{eq_t2}, \eqref{eq_INUW}
	are known in the literature as H\"older-logarithmic energy and regularity dependent stability estimates. 
			For the case when $\tau =1$
	in \eqref{eq_Theorem1}, \eqref{eq_Theorem2} or when $s=0$ in \eqref{eq_t2} the history 
	of such estimates in dimension $d=3$ goes back to \cite{Novikov2005+}, \cite{Novikov2008}, where 
	such energy and regularity dependent rapidly convergent approximate stability estimates were given for the inverse scattering problem.
	
	Then for Problem 1.1 energy dependent stability estimates changing from logarithmic type to H\"older type for high energies 
	were given in \cite{Isakov2011}. However, this high energy stability increasing of  \cite{Isakov2011} is slow.
		The studies of \cite{Novikov2005+}, \cite{Novikov2008}, \cite{Novikov2011}, \cite{Isakov2011} were continued, in particular, in \cite{NUW2011}, 
	\cite{IN2012}, \cite{INUW2014} and in the present work.
\end{Remark}

\begin{Remark}
In Theorems \ref{Theorem1}, \ref{Theorem2} we consider the  three-dimensional case  for simplicity only. Similar results
hold in dimension $d>3$.

		As regards to logarithmic and H\"older-logarithmic stability estimates for Problem 1.1 in dimension $d=2$, 
		we refer to 		\cite{NS2010}, \cite{Santacesaria2012}, \cite{Santacesaria2013}.		
		In addition, for problems like Problem 1.1 the history of energy and regularity dependent rapidly convergent approximate stability estimates in dimension $d=2$
		 goes back to \cite{Novikov1998}.
\end{Remark}


\begin{Remark}
In a similar way with results of \cite{IN2013}, \cite{IN2013+} and subsequent studies of \cite{PZ2013}, 
estimates 
\eqref{eq_Theorem1}, \eqref{eq_Theorem2} can be extended to the case when we do not assume
that condition (\ref{correct}) is fulfiled and consider an appropriate impedance boundary map (Robin-to-Robin map)
instead of the Dirichlet-to-Neumann map.
\end{Remark}

\begin{Remark}
  Apparently,  estimates analogous to estimates of Theorems \ref{Theorem1} and \ref{Theorem2}
  hold if we replace the difference of DtN maps by the difference of corresponding near field scattering data
  in a similar way with results of \cite{HH2001}, \cite{Isaev2013+},  \cite{IN2013++}. 
\end{Remark}


\begin{Remark}
	The optimality (in different senses)  of estimates like \eqref{eq_t1}, \eqref{eq_t2} was proved 
	in \cite{Mandache2001}, \cite{Isaev2011}, \cite{Isaev2013}. See also 
	\cite{CR2003}, \cite{Isaev2013++} and references therein for the case of inverse scattering problems.
\end{Remark}

\begin{Remark}
Estimates \eqref{eq_Theorem1}, \eqref{eq_Theorem2} for $\tau=1$
are roughly speaking coherent with stability properties of the approximate 
monochromatic inverse scattering reconstruction of \cite{Novikov2005+}, \cite{Novikov2008}, 
 implemented numerically in \cite{ABR2008}.
Estimates \eqref{eq_Theorem1}, \eqref{eq_Theorem2} for $E=0$
are roughly speaking coherent with stability properties of the reconstruction of \cite{Novikov2009}.

In addition, estimates \eqref{eq_Theorem1}, \eqref{eq_Theorem2}
can be used for the convergence rate analysis for iterative regularized reconstructions 
for Problem 1.1 in the framework of an effectivization of the approach 
of \cite{HH2001} for monochromatic inverse scattering problems.
\end{Remark}


\section{Lemmas}

Let $\hat{v}$ denote the Fourier transform of $v$:
\begin{equation}\label{eq_four}
	\hat{v} (\xi) = \mathcal{F} v (\xi)= \frac{1}{(2\pi)^3}\int\limits_{\mathbb{R}^3} e^{i\xi x} v(x) dx,
	\ \ \ \xi \in \mathbb{R}^3.
\end{equation}

\begin{Lemma}\label{Lemma_2.1}
	Suppose that $D$ satisfies (\ref{eq_c}) and  $v_1, v_2$ satisfy (\ref{eq_c1}), (\ref{correct})
	for some real $E$. 	Suppose also that $||v_{j}||_{L^\infty(D)} \leq N,\  j = 1,2$, for some $N > 0$. 
	Let $\delta$ be defined by (\ref{Theorem_delta}). Then 
			\begin{equation}\label{eq_Lemma2.1}
						|\hat{v}_2(\xi)-\hat{v}_1(\xi)| \leq c_1 (1+N)^2 \left(e^{2\rho L} \delta
					+ \frac{\|v_1-v_2\|_{L^{2}(D)}}{\sqrt{E+\rho^2}}\right)
			\end{equation}
	for any $\rho >0$ such that
	\begin{equation*}
								|\xi|\leq 2\sqrt{E+\rho^2},  \ \ \ E+\rho^2 \geq (1+N)^2 r_1^2, 
	\end{equation*}
where  $L = \max\limits_{x\in \partial D} |x|$ and
constants $c_1, r_1>0$ depend on $D$ only.
\end{Lemma}
Some version of estimate \eqref{eq_Lemma2.1} was given in \cite{IN2012} (see formula (4.13) of \cite{IN2012}).
Lemma \ref{Lemma_2.1} is proved in Section \ref{S_Proofs}. 
This proof is based on results presented in Section \ref{S_Faddeev}.

\begin{Lemma}\label{Lemma_2.2}
	Let  $w \in H^m(\mathbb{R}^3)$, $\|w\|_{H^m(\mathbb{R}^3)} 
\leq N_{H^m}$ for some real $m>0$ and $N_{H^m}>0$, where the space $H^m(\mathbb{R}^3)$ is defined in (\ref{def_Hm}). Then, for any $r>0$, 
				\begin{equation}\label{eq_Lemma2.2A1}
						\left(\int\limits_{|\xi|\geq r} |\mathcal{F}w(\xi)|^2 d\xi\right)^{1/2} \leq c_2 N_{H^m}\, r^{-m},
			\end{equation}
			 where 
			 $\mathcal{F} w$ is defined according to (\ref{eq_four}) and $c_2 = (2\pi)^{-3/2}$.
\end{Lemma}

\begin{Proof} {\it Lemma \ref{Lemma_2.2}.} 
Note that
	\begin{equation}\label{L1}
		\int\limits_{|\xi|\geq r} |\mathcal{F} w(\xi)|^2 d\xi
		\leq  \left\|\frac{(1+|\xi|^2)^{\frac{m}{2}} \mathcal{F} w }{r^m} \right\|_{L^2(\mathbb{R}^3)}^2.
	\end{equation}
		Using  (\ref{def_Hm}), \eqref{L1} and the Parseval theorem 
		\begin{equation}\label{Parseval}
		\|\mathcal{F}\tilde{w} \|_{L^2(\mathbb{R}^3)}  = (2\pi)^{-3/2}\|\tilde{w}\|_{L^2(\mathbb{R}^3)} 
		\end{equation}
			for $\tilde{w}  \equiv \mathcal{F}^{-1}(1+|\xi|^2)^{\frac{m}{2}} \mathcal{F} w$,
		 we get estimate \eqref{eq_Lemma2.2A1}.

\end{Proof}


\begin{Lemma}\label{Lemma_2.3}
	Let  $w \in W^m(\mathbb{R}^3)$, $\|w\|_{W^m(\mathbb{R}^3)} 
\leq N_{W^m}$ for some real $m>3$ and $N_{W^m}>0$, where the space $W^m(\mathbb{R}^3)$ is defined in (\ref{def_Wm}). Then, for any $r>0$, 
			\begin{equation}\label{eq_Lemma2.2A2}
							\int\limits_{|\xi|\geq r} |\mathcal{F} w(\xi)| d\xi \leq \tilde{c}_2 \frac{N_{W^m}}{m-3}\, r^{3-m},
			\end{equation}
				where 
			 $\mathcal{F} w$ is defined according to (\ref{eq_four}) and $\tilde{c}_2 = 4\pi$.
\end{Lemma}

\begin{Proof} {\it Lemma \ref{Lemma_2.3}.} 
Note that 
	\begin{equation}\label{L2}
			r^m |\mathcal{F} w(\xi)| \leq  (1+|\xi|^2)^{m/2} |\mathcal{F} w(\xi)| \leq   N_{W^m}
			 \ \ \text{ 	for  $|\xi| \geq r$.}
	\end{equation}
	Using \eqref{L2},  we obtain that
			\begin{equation}
			\int\limits_{|\xi|\geq r} |\mathcal{F}w(\xi)| d\xi  
			 \leq \int\limits_{r}\limits^{+\infty} \frac{N_{W^m}}{ t^m} \,  {4\pi t^2} dt
				\leq \frac{4\pi N_{W^m}}{m-3}\, r^{3-m}.
		\end{equation}

\end{Proof}

\section{ Proof of Theorem \ref{Theorem1}}

Using the Parseval formula \eqref{Parseval}, we get that 
\begin{equation}\label{L2_v2-v1}
	\|v_2-v_1\|_{L^{2}(D)} = (2\pi)^{3/2}\|\hat{v}_2-\hat{v}_1\|_{L^{2}(\mathbb{R}^3)} \leq (2\pi)^{3/2} (I_1(r) + I_2(r)),
\end{equation}
for $r>0$, where $\hat{v}_j$ is defined according to \eqref{eq_four} with $v_j\equiv 0$ on $\mathbb{R}^3 \setminus D$, $j=1,2$,
\begin{equation*}
	I_1(r) =	\left(\int\limits_{|\xi| \leq r} |\hat{v}_2(\xi)-\hat{v}_1(\xi)|^2 d\xi\right)^{1/2},
\end{equation*}
\begin{equation*}
	I_2(r) =	\left(\int\limits_{|\xi| \geq r} |\hat{v}_2(\xi)-\hat{v}_1(\xi)|^2 d\xi\right)^{1/2}.
\end{equation*}

Let 
\begin{equation}\label{L2_r}
	r  = q (1+N)^{-4/3}(E+\rho^2)^{1/3} , \ \ \ \  q = \frac{1}{2\pi}\left(\frac{16 \pi c_1^2}{3}\right)^{-1/3}, 
\end{equation}
where $c_1$ is the constant of Lemma \ref{Lemma_2.1}.

Then, using Lemma \ref{Lemma_2.1} for $|\xi| \leq r$, 
we get that
\begin{equation}\label{L2_I1}
	\begin{aligned}
		I_1(r) \leq  \left(\frac{4\pi r^3}{3}  c_1^2 (1+N)^4 \left(e^{2\rho L} \delta
					+ \frac{\|v_1-v_2\|_{L^{2}(D)}}{\sqrt{E+\rho^2}}\right)^2 \right)^{1/2} \leq\\
					\leq (2\pi)^{-3/2} \left( \frac{\sqrt{E+\rho^2}\, e^{2\rho L}\delta}{2}  + \frac{\|v_1-v_2\|_{L^{2}(D)}}{2}\right)
	\end{aligned}
\end{equation}
 for  $q (1+N)^{-4/3} (E+\rho^2)^{1/3}  \leq 2 \sqrt{E+\rho^2}$  and  $E+\rho^2 \geq (1+N)^2 r_1^2$.

In addition, using  \eqref{eq_Lemma2.2A1},  we have that
\begin{equation}\label{L2_I2}
	I_2(r) \leq c_2 N_{H^m}\, r^{-m}. 
\end{equation}
 Let $r_2=r_2(D)\geq r_1$ be such that
\begin{equation}\label{new_3.5}	
	E+\rho^2 \geq r_2^2 \ \ \Longrightarrow 	\ \ 
	\begin{aligned}
		q (E+\rho^2)^{1/3}  \leq 2 \sqrt{E+\rho^2}.
	\end{aligned}
\end{equation} 
Using \eqref{L2_v2-v1}, \eqref{L2_I1}--\eqref{new_3.5} with $r$ defined in \eqref{L2_r}, we obtain that 
	\begin{equation}
		\begin{aligned}
		\|v_2-v_1\|_{L^{2}(D)} \leq \frac{ \sqrt{E+\rho^2}\,  e^{2\rho L} \delta }{2}   + \frac{\|v_1-v_2\|_{L^{2}(D)}}{2}+ \\
		+ (2\pi)^{3/2} c_2  \frac{(1+N)^{\frac{4m}{3}}}{q^{m}} N_{H^m} (E+\rho^2)^{-\frac{m}{3}}, 
		\end{aligned}
	\end{equation}
	\begin{equation}\label{L2_rho}
		\begin{aligned}
		\frac{1}{2}\|v_2-v_1\|_{L^{2}(D)} \leq &\frac{\sqrt{E+\rho^2}\,  e^{2\rho L} \delta}{2}  + \\
		&+  \frac{(1+N)^{\frac{4m}{3}}}{q^{m}} N_{H^m} (E+\rho^2)^{-\frac{m}{3}}, 
		\end{aligned}
	\end{equation}
	for $E+\rho^2 \geq  (1+N)^2 r_2^2$, where $L,c_2$ are the constants of Lemmas \ref{Lemma_2.1}, \ref{Lemma_2.2} and $q, r_2$ are the constants of formulas \eqref{L2_r}, \eqref{new_3.5}.

Let $\tau \in (0,1)$ and 
\begin{equation}\label{L2_gamma}
	\gamma = \frac{1-\tau}{2L}, \ \ \ \  \rho = \gamma\ln\left(3+ \delta^{-1}\right). 
\end{equation} 
 Due to (\ref{L2_rho}), for $\delta$  such that
\begin{equation}\label{fa_1}
	 E+ \left(\gamma\ln (3+ \delta^{-1})\right)^2  \geq 	(1+N)^2 r_2^2, 
\end{equation} 
   the following estimate holds:
\begin{equation}\label{L2_est}
	\begin{aligned}
		\frac{1}{2}\|v_1 - v_2&\|_{L^{2}(D)} \leq \\ &\leq 
		\frac{1}{2}\left(E+ \left(\gamma\ln\left(3+ \delta^{-1}\right)\right)^2\right)^{1/2}
		\left(3+ \delta^{-1}\right)^{2\gamma L} \delta + \\
		&+    \frac{(1+N)^{\frac{4m}{3}}}{q^{m}} N_{H^m}\left(E+\left(\gamma\ln\left(3+ \delta^{-1}\right)\right)^2 \right)^{-\frac{m}{3}},  
	\end{aligned}
\end{equation}
where  $\gamma$ is defined  in (\ref{L2_gamma}). 
Note that 
\begin{equation}\label{fa_2}
	\left(3+ \delta^{-1}\right)^{2\gamma L} \delta = \left(1+ 3\delta \right)^{1-\tau} \delta^{\tau} \leq 4 \delta^{\tau} 
	\ \ \ \text{ for $\delta\leq 1$.}	
\end{equation}
Combining (\ref{L2_est}), \eqref{fa_2}, we get that
\begin{equation}\label{L2_almost1}
	\begin{aligned}
	||v_2 - v_1||_{L^2 (D)} \leq  
	A_1\left( \lambda\left( E + \gamma^2 \left(\ln\left(3+\delta ^{-1}\right)\right)^2 \right) \right)^{\frac{1}{2}} \delta^\tau + \\ 
	+
	 B_1  \left( 1+N\right)^{\frac{4m}{3}} N_{H^m}  \left( \lambda\left( E + \gamma^2 \left(\ln\left(3+\delta ^{-1}\right)\right)^2 \right) \right)^{-\frac{m}{3}} 
	\end{aligned}
\end{equation}	
	for $\delta \leq 1$ satisfying \eqref{fa_1} and some positive constants $A_1,B_1,\lambda$ depending on $D$ only.

In view of definition \eqref{def_Hm}, we have that  
\begin{equation*}
||v_2 - v_1||_{L^2 (D)} \leq ||v_2 - v_1||_{H^m (\mathbb{R}^3)} \leq N_{H^m}.
\end{equation*}
Hence, we get that, for  $0 < E+ \left(\gamma\ln (3+ \delta^{-1})\right)^2  \leq 	(1+N)^2 r_2^2$,
\begin{equation}\label{L2_almost2}
	||v_2 - v_1||_{L^2 (D)}  \leq \left( 1+N\right)^{\frac{4m}{3}} N_{H^m}  \left( \frac{E + \gamma^2 \left(\ln\left(3+\delta ^{-1}\right)\right)^2}{r_2^2} \right)^{-\frac{m}{3}}. 
\end{equation}
On other hand, in the case when 
$E+ \left(\gamma\ln (3+\delta^{-1})\right)^2  \geq 	(1+N)^2 r_2^2$ and $\delta > 1$ we have that
\begin{equation}\label{L2_almost3}
	\begin{aligned}
	||v_2 - v_1||_{L^2 (D)}  \leq c_3 ||v_2 - v_1||_{L^\infty (D)} \leq c_3 2N \leq \\
	\leq 2c_3 \left( \frac{E + \gamma^2 \left(\ln\left(3+\delta ^{-1}\right)\right)^2}{r_2^2} \right)^{\frac{1}{2}} \delta^\tau, 
	\end{aligned}
\end{equation}
where 
\begin{equation}\label{new_3.14}
c_3 = \left( \int\limits_D 1\, dx \right)^{1/2}.
\end{equation}

Combining \eqref{L2_gamma}, \eqref{L2_almost1}--\eqref{L2_almost3}, we obtain estimate \eqref{eq_Theorem1}. This completes the proof of Theorem \ref{Theorem1}.  


\section{Proof of Theorem \ref{Theorem2} }

Due to the inverse Fourier transform formula
\begin{equation}
	v(x)  = \int\limits_{\mathbb{R}^3} e^{-i\xi x}\hat{v}(\xi) d\xi, \ \ \ x\in \mathbb{R}^3, 
\end{equation}
we have that 
\begin{equation}\label{L0_v2-v1}
	\begin{aligned}
		\|v_1 - v_2\|_{L^{\infty}(D)} \leq 
		\sup\limits_{x\in D}\left|\int\limits_{\mathbb{R}^3} 
		e^{-i\xi x}\left(\hat{v}_2(\xi) - \hat{v}_1(\xi)\right) d\xi\right| 
		\leq \tilde{I}_1(r) + \tilde{I}_2(r) 
	\end{aligned}
\end{equation}
for $r>0$, where
\begin{equation*}
	\begin{aligned}
		\tilde{I}_1(r) = \int\limits_{|\xi|\leq r} |\hat{v}_2(\xi) - \hat{v}_1(\xi)| d\xi, \\
		\tilde{I}_2(r) = \int\limits_{|\xi|\geq r} |\hat{v}_2(\xi) - \hat{v}_1(\xi)| d\xi.
	\end{aligned}
\end{equation*}

Let \begin{equation}\label{L0_r}
	r  = \tilde{q} (1+N)^{-2/3}(E+\rho^2)^{1/6} , \ \ \ \  
	\tilde{q} = \left(\frac{8 \pi c_1 c_3}{3}\right)^{-1/3}, 
\end{equation}
where $c_1$ is the constant of Lemma \ref{Lemma_2.1} and $c_3$ is defined by \eqref{new_3.14}.

Then, combining the definition of $\tilde{I}_1$, Lemma \ref{Lemma_2.1} for $|\xi| \leq r$ and 
the inequality $$||v_2 - v_1||_{L^2 (D)}  \leq c_3 ||v_2 - v_1||_{L^\infty (D)},$$ 
we get that
\begin{equation}\label{L0_I1}
	\begin{aligned}
		\tilde{I}_1(r) \leq  \frac{4\pi r^3}{3}  c_1 (1+N)^2 \left(e^{2\rho L} \delta
					+ \frac{c_3 \|v_1-v_2\|_{L^{\infty}(D)}}{\sqrt{E+\rho^2}}\right)  \leq\\
					 \leq \frac{1}{2c_3}  \sqrt{E+\rho^2}\, e^{2\rho L}\delta + \frac{\|v_1-v_2\|_{L^{\infty}(D)}}{2}
	\end{aligned}
\end{equation}
 for  $\tilde{q} (1+N)^{-2/3} (E+\rho^2)^{1/6}  \leq 2 \sqrt{E+\rho^2}$  and  $E+\rho^2 \geq (1+N)^2 r_1^2$.

In addition, using  \eqref{eq_Lemma2.2A2},  we get that
\begin{equation}\label{L0_I2}
	\tilde{I}_2(r) \leq \tilde{c}_2 \frac{N_{W^m}}{m-3}\, r^{3-m}. 
\end{equation}

 Let $\tilde{r}_2=\tilde{r}_2(D)\geq r_1$ be such that
\begin{equation}\label{new_4.5}	
	E+\rho^2 \geq \tilde{r}_2^2 \ \ \Longrightarrow 	\ \ 
	\begin{aligned}
		\tilde{q} (E+\rho^2)^{1/6}  \leq 2 \sqrt{E+\rho^2}.
	\end{aligned}
\end{equation} 
Using \eqref{L0_v2-v1}, \eqref{L0_I1}--\eqref{new_4.5} with $r$ defined in \eqref{L0_r}, we obtain that 
	\begin{equation}
		\begin{aligned}
		\|v_2-v_1\|_{L^{\infty}(D)} \leq \frac{1}{2c_3}
		 \sqrt{E+\rho^2}\,  e^{2\rho L} \delta  + \frac{\|v_1-v_2\|_{L^{\infty}(D)}}{2}+ \\
		+  \tilde{c}_2  \frac{(1+N)^{\frac{2(m-3)}{3}}}{(m-3)\tilde{q}^{m-3}} N_{W^m} (E+\rho^2)^{-\frac{m-3}{6}},
		\end{aligned}
	\end{equation}
	\begin{equation}\label{L0_rho}
		\begin{aligned}
		\frac{1}{2}\|v_2-v_1\|_{L^{\infty}(D)} &\leq \frac{1}{2c_3}
		 \sqrt{E+\rho^2}\,  e^{2\rho L} \delta  + \\
		&+  4\pi  \frac{(1+N)^{\frac{2(m-3)}{3}}}{(m-3)\tilde{q}^{m-3}} N_{W^m} (E+\rho^2)^{-\frac{m-3}{6}}
		\end{aligned}
	\end{equation}
	for $E+\rho^2 \geq  (1+N)^2 \tilde{r}_2^2$, where  $L$, $\tilde{c_2}$ 
	are the constants of Lemmas \ref{Lemma_2.1},
	\ref{Lemma_2.3} and $c_3$, $\tilde{q}$, $\tilde{r}_2$ are the constants of formulas \eqref{new_3.14}, 
	\eqref{L0_r}, \eqref{new_4.5}.
	
Let $\tau \in (0,1)$ and 
\begin{equation}\label{L0_gamma}
	\gamma = \frac{1-\tau}{2L}, \ \ \ \  \rho = \gamma\ln\left(3+ \delta^{-1}\right). 
\end{equation} 
  Due to (\ref{L0_rho}), for  $\delta$ such that
\begin{equation}\label{0fa_1}
	 E+ \left(\gamma\ln (3+\delta^{-1})\right)^2  \geq 	(1+N)^2 \tilde{r}_2^2, 
\end{equation} 
   the following estimate holds:
\begin{equation}\label{L0_est}
	\begin{aligned}
		\frac{1}{2}\|v_1 - v_2\|_{L^{\infty}(D)} &\leq \\ \leq 
		\frac{1}{2c_3}&\left(E+ \left(\gamma\ln\left(3+ \delta^{-1}\right)\right)^2\right)^{1/2}
		\left(3+ \delta^{-1}\right)^{2\gamma L} \delta + \\
		+   &4\pi  \frac{(1+N)^{\frac{2(m-3)}{3}}}{(m-3)\tilde{q}^{m}} N_{W^m}\left(E+\left(\gamma\ln\left(3+ \delta^{-1}\right)\right)^2 \right)^{-\frac{m-3}{6}},  
	\end{aligned}
\end{equation}
where  $\gamma$ is defined  in (\ref{L0_gamma}). 
Note that 
\begin{equation}\label{0fa_2}
	\left(3+ \delta^{-1}\right)^{2\gamma L} \delta = \left(1+ 3\delta \right)^{1-\tau} \delta^{\tau} \leq 4 \delta^{\tau}
	\ \ \text{ for $\delta\leq 1$.}
\end{equation}
Combining (\ref{L0_est}), \eqref{0fa_2}, we get that
\begin{equation}\label{L0_almost1}
	\begin{aligned}
	||v_2 - v_1||_{L^\infty (D)} \leq  
	\tilde{A}_1\left( \tilde{\lambda}\left( E + \gamma^2 \left(\ln\left(3+\delta ^{-1}\right)\right)^2 \right) \right)^{\frac{1}{2}} \delta^\tau + \\ 
	+
	 \tilde{B}_1  \frac{ \left( 1+N\right)^{\frac{2(m-3)}{3}} N_{W^m}}{m-3}  \left( \tilde{\lambda}\left( E + \gamma^2 \left(\ln\left(3+\delta ^{-1}\right)\right)^2 \right) \right)^{-\frac{m-3}{6}} 
	\end{aligned}
\end{equation}	
	for $\delta \leq 1$ satisfying \eqref{0fa_1} and some positive constants $\tilde{A}_1,\tilde{B}_1,\tilde{\lambda}$ depending on $D$ only.	
	
	Using \eqref{def_Wm} and \eqref{L0_v2-v1}, we get that 
\begin{equation}\label{new_4.13}
	\begin{aligned}
||v_2 - v_1||_{L^\infty (D)} \leq 
    \int\limits_{\mathbb{R}^3} \left( {(1+|\xi|^2)^{-m/2}} ||v_2 - v_1||_{W^m(\mathbb{R}^3)} \right) d\xi \leq\\
    \leq   N_{W^m} \int\limits_0\limits^{+\infty} \frac{4\pi t^2}{(1+t^2)^{m/2}} dt  
    \leq  c_4 \frac{e^{m-3}}{m-3} N_{W^m}
\end{aligned}
\end{equation}
for some $c_4>0$. Here we used also that
\begin{equation*}
	\int\limits_0\limits^{+\infty} \frac{4\pi t^2}{(1+t^2)^{m/2}} dt \leq 
	\int\limits_0\limits^{1} {4\pi t^2} dt + 
	\int\limits_1\limits^{+\infty} \frac{4\pi t^2}{t^m} dt \leq c_4 \left(1+ \frac{1}{m-3}\right) \leq  
	c_4 \frac{e^{m-3}}{m-3} .
\end{equation*}
Using \eqref{new_4.13}, we get that, for  $0 < E+ \left(\gamma\ln (3+\delta^{-1})\right)^2  \leq 	(1+N)^2 \tilde{r}_2^2$,
\begin{equation}\label{L0_almost2}
	\begin{aligned}
	||v_2 - v_1||_{L^\infty (D)}  &\leq \\
	\leq c_4 & \frac{\left( 1+N\right)^{\frac{2(m-3)}{3}} N_{W^m}}{m-3}  \left( \frac{E + \gamma^2 \left(\ln\left(3+\delta ^{-1}\right)\right)^2}{e^6 \tilde{r}_2^2}  \right)^{-\frac{m-3}{6}}. 
	\end{aligned}
\end{equation}
On other hand, in the case when 
$E+ \left(\gamma\ln \delta^{-1}\right)^2  \geq 	(1+N)^2 \tilde{r}_2^2$ and $\delta > 1$ we have that
\begin{equation}\label{L0_almost3}
	\begin{aligned}
	||v_2 - v_1||_{L^\infty (D)}  \leq  2N 
	\leq  2\left( \frac{E + \gamma^2 \left(\ln\left(3+\delta ^{-1}\right)\right)^2}{\tilde{r}_2^2} \right)^{\frac{1}{2}} \delta^\tau. 
	\end{aligned}
\end{equation}

Combining \eqref{L0_gamma}, \eqref{L0_almost1}, \eqref{L0_almost2} and \eqref{L0_almost3}, we obtain estimate \eqref{eq_Theorem2}. This completes the proof of Theorem \ref{Theorem2}.  

\section{Faddeev functions}\label{S_Faddeev}
Suppose that 
\begin{equation}\label{O4.7}
	v \in L^{\infty}(D), \ \ \ v \equiv 0 \text{ on } \mathbb{R}^3\setminus D,
\end{equation}
where $D$ satisfies \eqref{eq_c}. More generally, one can assume that
\begin{equation}\label{O4.6}
	\begin{aligned}
	v \text{ is a sufficiently regular function on } \mathbb{R}^3 \\
	\text{ with sufficient decay at infinity.}
	\end{aligned}
\end{equation}

Under assumptions \eqref{O4.6}, we consider the  functions $\psi$, $\mu$, $h$: 
\begin{equation}\label{O4.1}
	\psi(x,k) = e^{ikx} \mu(x,k),
\end{equation}
\begin{equation}\label{O4.2}
	\begin{aligned}
	\mu(x,k) = 1 + \int\limits_{\mathbb{R}^3} g(x-y,k) v(y)\mu(y,k) dy, \\
  g(x,k) = - (2\pi)^{-3} \int\limits_{\mathbb{R}^3} \frac{e^{i\xi x} d\xi}{\xi^2 + 2k\xi},
	\end{aligned}
\end{equation}
where $x\in \mathbb{R}^3$, $k\in \mathbb{C}^3$, $\mbox{Im}\, k\neq 0$, 
\begin{equation}\label{O4.4}
	h(k,l) = (2\pi)^{-3} \int\limits_{\mathbb{R}^3} e^{i(k-l)x}v(x) \mu(x,k) dx,
\end{equation}
where 
$
	 k,l \in \mathbb{C}^3, \ k^2=l^2,\ \mbox{Im}\, k = \mbox{Im}\, l \neq 0.
$
Here,  \eqref{O4.2} at fixed $k$ is considered as a linear integral equation for $\mu$,  where
$\mu$ is sought in $L^{\infty}(\mathbb{R}^3)$.

The functions  $\psi$, $h$  and $G = e^{ikx}g$ are known as the Faddeev functions, see \cite{Faddeev1965}, \cite{Faddeev1974}, \cite{Henkin1987}, \cite{Novikov 1988}.
These functions were introduced for the first time  in \cite{Faddeev1965}, \cite{Faddeev1974}.

In particular, we have that
\begin{equation*}
	(\Delta+k^2) G(x,k) = \delta(x), 
\end{equation*}
\begin{equation*}
	(-\Delta+v(x)) \psi(x,k) = k^2\psi(x,k), 
\end{equation*}
where $x\in\mathbb{R}^3$,  $k \in \mathbb{C}^3\setminus \mathbb{R}^3$.

We recall also that   the Faddeev functions $G$, $\psi$, $h$ are some extension to the complex domain of functions     of the classical scattering theory for the Schr\"odinger equation (in particular, $h$ is an extension of the classical scattering amplitude).

 Note also that $G$, $\psi$, $h$ in their zero energy restriction, that is for $k^2=0$, $l^2=0$,
  were considered for the first time in \cite{Beals1985}.
The Faddeev functions $G$, $\psi$, $h$ were, actually, rediscovered in \cite{Beals1985}.

For further considerations we will use the following notations: 
\begin{equation*}
\begin{aligned}
		\Sigma_E = \left\{ k\in \mathbb{C}^3: k^2 = k_1^2 + k_2^2 + k_3 ^2 = E\right\},\\
		\Theta_E = \left\{ k\in \Sigma_E,\  l\in\Sigma_E: \mbox{Im}\, k = \mbox{Im}\, l\right\},\\
		|k| = (|\mbox{Re}\,k|^2 +|\mbox{Im}\,k|^2)^{1/2} \  \text{ for } k \in \mathbb{C}^3.
\end{aligned}		
\end{equation*}
Under assumptions \eqref{O4.6}, we have that:
\begin{equation}\label{mu_1}
	\mu(x,k) \rightarrow 1 \ \  \text{ as } \ \ |k|\rightarrow \infty, 
\end{equation}
where $x\in \mathbb{R}^3$, $k \in \Sigma_E$;
\begin{equation}\label{lim_1}
			\hat{v}(\xi) = \lim\limits_
		{\scriptsize
			\begin{array}{c}
			(k,l)\in \Theta_E,\, k-l=\xi\\
			|\mbox{Im}\,k|=|\mbox{Im}\,l|\rightarrow \infty
			\end{array}
		} h(k,l)\ \ \ 	 \text{ for any } \xi \in \mathbb{R}^3,
	\end{equation}
		where $\hat{v}$ is defined by \eqref{eq_four}.

Results of the type (\ref{mu_1}) go back to \cite{Beals1985}.  
Results of the type (\ref{lim_1}) go back to \cite{Henkin1987}. 
These results follow, for example, from 
equation (\ref{O4.2}), formula (\ref{O4.4}) and 
the following estimates:
\begin{equation}\label{3.11}
	\begin{aligned}
 g(x, k) = O(|x|^{-1}) 
 \  \text{ for }  x \in \mathbb{R}^3,\\
		\text{ uniformly in } k\in \mathbb{C}^3\setminus \mathbb{R}^3,
	\end{aligned}
\end{equation}

\begin{equation}\label{4.15.new}
	\begin{aligned}
	\| \Lambda^{-s} g(k) \Lambda^{-s}\|_{L^2(\mathbb{R}^3)\rightarrow L^2(\mathbb{R}^3)} = O(|k|^{-1}),
	\text{ for } s>1/2,
	 \\ \text{ as } \ |k|\rightarrow \infty,\ \ \
		k\in \mathbb{C}^3\setminus \mathbb{R}^3,
	\end{aligned}
\end{equation}
where $g(x,k)$ is defined in \eqref{O4.2}, $g(k)$ denotes the integral operator with the Schwartz kernel $g(x-y,k)$ 
 and $\Lambda$ denotes the multiplication operator by the function $(1+|x|^2)^{1/2}$. 
Estimate  \eqref{3.11} was given in \cite{Henkin1987}.
Estimate (\ref{4.15.new}) was formulated, first, in \cite{LN1987}. Concerning proof of (\ref{4.15.new}), see \cite{Weder1991}.  In addition, 
estimate (\ref{4.15.new}) in its zero energy restriction goes back to \cite{SU1987}.

In the present work we use the following lemma:
\begin{Lemma}\label{Lemma_3.1}
 Let $D$ satisfy (\ref{eq_c}) and $v$ satisfy (\ref{O4.7}).
Let $||v||_{L^\infty(D)} \leq N$ for some $N>0$. Then 
\begin{equation}
	 |\mu(x,k)|  \leq c_5(1+N) \ \ \ \text{ for } \ \ x \in \mathbb{R}^3, \ |k|\geq r_3(1+N),
\end{equation}
where  $ \mu(x,k)$ is the Faddeev function of (\ref{O4.2}) and 
constants $c_5, r_3>0$ depend on $D$ only.  
\end{Lemma}

Lemma \ref{Lemma_3.1} is proved in Section 6. This proof is based on  estimates \eqref{3.11} and (\ref{4.15.new}).

In addition, we have that (see \cite{Novikov1996}, \cite{Novikov2005}):
\begin{equation}\label{delta_h0}
	\begin{aligned}
	h_2(k,l) - h_1(k,l) = (2\pi)^{-3} \int\limits_{\mathbb{R}^3} \psi_1(x,-l) (v_2(x) - v_1(x)) \psi_2(x,k) dx 	\\ 	 
	\text{ for }  (k,l) \in \Theta_E,\ |\mbox{Im}\,k| = |\mbox{Im} \,l| \neq 0,\\ \text{ and $v_1$, $v_2$ satisfying (\ref{O4.6}),}
	\end{aligned}
\end{equation}
\begin{equation}\label{delta_h}
	\begin{aligned}
	h_2(k,l) - h_1(k,l) = (2\pi)^{-3} \int\limits_{\partial D} \psi_1(x,-l) 
	\left[\left(\hat{\Phi}_{2} - \hat{\Phi}_{1}\right) 
	\psi_2(\cdot,k)\right](x) dx\\
	\text{ for }  (k,l) \in \Theta_E,\ |\mbox{Im}\,k| = |\mbox{Im} \,l| \neq 0,\\ \text{ and $v_1$, $v_2$ satisfying (\ref{correct}), (\ref{O4.7}),}
	\end{aligned}
\end{equation}
where  $\psi_j$, $h_j$ denote $\psi$ and $h$  of  (\ref{O4.1}) and (\ref{O4.4}) for $v = v_j$, and $\hat{\Phi}_{j}$ denotes
the Dirichlet-to-Neumann map $\hat{\Phi}$ for $v = v_j$ in $D$, where $j=1,2$.

In the present work we also use the following lemma:
\begin{Lemma}\label{Lemma_3.2}
	Let $D$ satisfy (\ref{eq_c}). Let $v_j$ satisfy (\ref{O4.7}), 
	$||v_{j}||_{L^\infty(D)} \leq N,\  j = 1,2$, for some $N > 0$. 
 Then 
		\begin{equation}\label{eq_Lemma3.2}
	\begin{aligned}
		|\hat{v}_1(\xi) - \hat{v}_2(\xi) - h_1(k,l) + h_2(k,l)|
		\leq \frac{
		c_6 N(1+N)
		\|v_1-v_2\|_{L^{2}(D)}
		}
		{(E+\rho^2)^{1/2}}  	\\ \text{ for }  (k,l) \in \Theta_E,  \ \    \xi = k-l, \ \ 
		|{\rm Im}\,k| = |{\rm Im} \,l| = \rho,\\ E+\rho^2 \geq r_4^2(1+N)^2,
	\end{aligned}
	\end{equation}
where $E \in \mathbb{R}$, 
$\hat{v}_j$ is the Fourier transform of $v_j$, $h_j$ denotes  $h$  of  (\ref{O4.4})  for $v = v_j$, $(j=1,2)$ 
and constants $c_6, r_4>0$ depend on $D$ only.

\end{Lemma}

Some versions of estimate \eqref{eq_Lemma3.2} were given in 
\cite{IN2012}, \cite{Novikov1996}, \cite{Novikov2005} (see, for example, formula (3.18) of \cite{IN2012}).
 Lemma \ref{Lemma_3.2} is proved in Section 6.


\section{Proofs of Lemmas \ref{Lemma_2.1}, \ref{Lemma_3.1} and \ref{Lemma_3.2}}\label{S_Proofs}

\begin{Proof} {\it Lemma \ref{Lemma_3.1}.}
Using \eqref{O4.7}, \eqref{O4.2} and \eqref{4.15.new}, we get that 
\begin{equation}\label{4.2a}
	\begin{aligned}
	 \|\mu(\cdot,k) - 1\|_{L^2(D)} 
	 \leq \left\|\int\limits_{\mathbb{R}^3} g(\cdot-y) v(y) \mu (y,k)dy\right\|_{L^2(D)}\leq  \\
	 	\leq  c_7 \frac{N}{|k|} \|\mu(\cdot,k)\|_{L^2(D)}, 
	\end{aligned}
\end{equation}
\begin{equation}\label{6.2}
	\begin{aligned}
	 \|\mu(\cdot,k)\|_{L^2(D)} 
	 \leq c_3+  c_7 \frac{N}{|k|} \|\mu(\cdot,k)\|_{L^2(D)}, 
	\end{aligned}
\end{equation}
where $c_3$ is defined by \eqref{new_3.14} and $c_7$ is some positive 
constant depending on $D$ only.   Hence, we obtain that 
\begin{equation}\label{4.2}
		 \|\mu(\cdot,k)\|_{L^2(D)}  \leq 2c_3   \ \ \ \text{ for } \ \ \ |k| \geq 2c_7 N.
\end{equation}
We use also that 
\begin{equation}\label{6.4}
	\int\limits_D \frac{1}{|x-y|^2} dy \leq 
	\int\limits_D 1\, dy  + \int\limits_{|x-y| \leq 1} \frac{1}{|x-y|^2} dy \leq c_8^2, \ \ x\in D,
\end{equation}
where $c_8=c_8(D)>0$. Using \eqref{O4.7}, \eqref{O4.2}, \eqref{3.11}, \eqref{4.2}, \eqref{6.4},
 we get that


\begin{equation}\label{6.5}
	\begin{aligned}
	|\mu(x,k)| &\leq 1 + \left| \int\limits_{D} g(x-y) v(y) \mu (y,k)dy \right| \leq\\
	&\leq 1 +  \left(\int\limits_{D} |g(x-y)|^2 dy \right)^{1/2} N\,  \|\mu(\cdot,k)\|_{L^2(\mathbb{R}^3)} \leq\\
	 &\leq c_5(D)(1+N)  \ \hspace{10mm} \text{ for } \ \ x\in D, \ |k|\geq {2c_7N}. 
	\end{aligned}
\end{equation}
\end{Proof}

\begin{Proof} {\it Lemma \ref{Lemma_3.2}.}
	Due to \eqref{O4.7}, \eqref{delta_h0}, we have that
	\begin{equation}\label{4.4}
		\begin{aligned}
		h_2(k,l) - h_1(k,l)		
		= (2\pi)^{-3} \int\limits_{D} \psi_1(x,-l) (v_2(x) - v_1(x))  \psi_2(x,k) dx=
		\\=  (2\pi)^{-3} \int\limits_{D} e^{i(k-l)x} \mu_1(x,-l) (v_2(x) - v_1(x))  \mu_2(x,k) dx=
		\\ = \hat{v}_2(k-l) - \hat{v}_1(k-l) +  I_\Delta 
		\end{aligned}
	\end{equation}	 
 $\text{ for }  (k,l) \in \Theta_E,\ |\mbox{Im}\,k| = |\mbox{Im} \,l| \neq 0$, where 
	\begin{equation}\label{I_delta}
		\begin{aligned}
			I_\Delta = (2\pi)^{-3} \int\limits_{D} (\mu_1(x,-l) - 1)(v_2(x) - v_1(x)) \mu_2(x,k) dx + \\
			+  (2\pi)^{-3} \int\limits_{D} \mu_1(x,-l)(v_2(x) - v_1(x)) (\mu_2(x,k)-1) dx +\\
			+ (2\pi)^{-3} \int\limits_{D} (\mu_1(x,-l) - 1) (v_2(x) - v_1(x)) (\mu_2(x,k)-1) dx.
		\end{aligned}
	\end{equation}
	
	Note that, for $(k,l) \in \Theta_E$, $E\in \mathbb{R}$, $|\mbox{Im}\,k| = |\mbox{Im} \,l| = \rho$,
	\begin{equation}\label{4.6a}
		|k| = \sqrt{|\text{Re}\,k|^2 + |\text{Im}\,k|^2 } = \sqrt{k^2 + 2\, 
		|\text{Im}\,k|^2} = \sqrt{E+ 2\rho^2} =		 |l|.
	\end{equation}
	Using estimates \eqref{4.2a}, \eqref{4.2}, \eqref{6.5} in \eqref{I_delta}, we  get that 
	\begin{equation}\label{4.6}
		\begin{aligned}
	I_\Delta \leq   
	(2\pi)^{-3} \Bigg( \|\mu_1(\cdot,-l)-1\|_{L^2(D)} \|v_2-v_1\|_{L^2(D)} \|\mu_2(\cdot,-l)\|_{L^\infty(D)} +\\
	+  \|\mu_1(\cdot,-l)\|_{L^\infty(D)} \|v_2-v_1\|_{L^2(D)} \|\mu_2(\cdot,-l)-1\|_{L^2(D)} + \\
	+  \|\mu_1(\cdot,-l)-1\|_{L^\infty(D)} \|v_2-v_1\|_{L^2(D)} \|\mu_2(\cdot,-l)-1\|_{L^2(D)} \Bigg)\leq 
	\\
	\leq
	\frac{	2c_3c_7 N
		\|v_1-v_2\|_{L^{2}(D)}  c_5(1+N)	}
		{(2\pi)^3 |k|}  + 
		 \frac{c_5(1+N) \|v_1-v_2\|_{L^{2}(D)} 	2 c_3c_7 N	} 		{(2\pi)^3 |l|}  +\\
		+  \frac{	(1+ c_5(1+N))
		\|v_1-v_2\|_{L^{2}(D)}   2 c_3c_7 N
		}
		{(2\pi)^3|l|}  \leq \\  \leq c_8(D)  
		\frac{ N(1+N)		\|v_1-v_2\|_{L^{2}(D)}
		}
		{ \sqrt{E+ 2\rho^2}}
		\end{aligned}
	\end{equation}
		for $(k,l) \in \Theta_E,\ |\mbox{Im}\,k| = |\mbox{Im} \,l| = \rho$  and $|k|=|l|=\sqrt{E+2\rho^2} \geq 2c_7N$.
		
		Formula \eqref{4.4} and estimate \eqref{4.6} imply \eqref{eq_Lemma3.2}.
	
\end{Proof}

\begin{Proof} {\it Lemma \ref{Lemma_2.1}.}
Due to \eqref{delta_h0}, we have that
\begin{equation}\label{4.7}
	\begin{aligned}
		|h_2(k,l)- h_1(k,l)| \leq c_9 \|\psi_1(\cdot, -l)\|_{L^{\infty}(\partial D)}\, \delta \,  
		\|\psi_2(\cdot, k)\|_{L^{\infty}(\partial D)},\\
		 (k,l)\in \Theta_E,  \ \ |\text{Im}\,k|= |\text{Im}\,l| \neq 0,
	\end{aligned}
\end{equation}
where 
\begin{equation*}
	c_9 = (2\pi)^{-3} \int\limits_{\partial D} dx. 
\end{equation*}
Using formula \eqref{O4.1} and Lemma \ref{Lemma_3.1}, we find that
\begin{equation}\label{4.8}
	\begin{aligned}
	\|\psi_j(\cdot, k)\|_{L^{\infty}(\partial D)} \leq c_5(1+N) e^{|\text{Im}\,k| L},\  j=1,2, \\  \text{ for }\    k \in \Sigma_E, \  |k| \geq r_3(1+N),
	\end{aligned}
\end{equation}
where $L = \max\limits_{x\in \partial D} |x|$.
Combining \eqref{4.6a}, \eqref{4.7} and \eqref{4.8}, we get that
\begin{equation}\label{4.10}
	\begin{aligned}
	|h_2(k,l)- h_1(k,l)| \leq c_9 c_5^2 (1+N)^2 e^{2 \rho L} \delta, \\
		\text{ for } \ (k,l)\in \Theta_E, \ \rho = |\text{Im}\,k|= |\text{Im}\,l|,\\
		 E + \rho^2 \geq  r_3^2(1+N)^2.
	\end{aligned}
\end{equation}

Note that for any $\xi \in \mathbb{R}^3$ satisfying $|\xi| \leq 2\sqrt{E+ \rho^2}$ (where $\rho>0$) 
there exist some pair $(k,l)\in \Theta_E$ such that $\xi = k-l$ and $|\text{Im}\,k|= |\text{Im}\,l| = \rho$. 
Therefore, estimates  \eqref{eq_Lemma3.2} and \eqref{4.10} imply \eqref{eq_Lemma2.1}.
\end{Proof}


\noindent
{ {\bf M.I. Isaev}\\
Centre de Math\'ematiques Appliqu\'ees, Ecole Polytechnique,

91128 Palaiseau, France\\
Moscow Institute of Physics and Technology,

141700 Dolgoprudny, Russia\\
e-mail: \tt{isaev.m.i@gmail.com}}\\

\noindent
{ {\bf R.G. Novikov}\\
Centre de Math\'ematiques Appliqu\'ees, Ecole Polytechnique,

91128 Palaiseau, France\\
Institute of Earthquake Prediction Theory and Math. Geophysics RAS,

117997 Moscow, Russia\\ 
e-mail: \tt{novikov@cmap.polytechnique.fr}}

\end{document}